\documentclass{article}
\def\be{\begin{equation}}
\def\ee{\end{equation}}
\newcommand{\ff}[1]{{\mbox{\boldmath $#1$}}}
\def\a{\alpha}

\def\lam{\lambda}
\def\x{\ff{x}}

\begin{document}
\title{Cuckoo Search: Recent Advances and Applications}
\author{Xin-She Yang \\
School of Science and Technology, \\ Middlesex University, The Burroughs, \\
London NW4 4BT, UK.\\
\and
Suash Deb \\
Cambridge Institute of Technology,  \\ 
Cambridge Village, Tatisilwai, Ranchi-835103, \\
Jharkhand, India. }

\date{}

\maketitle

\begin{abstract}
Cuckoo search (CS) is a relatively new algorithm, developed by Yang and Deb in 2009, and
CS is efficient in solving global optimization problems. In this paper, we review
the fundamental ideas of cuckoo search and the latest developments as well
as its applications. We analyze the algorithm and gain insight into its search
mechanisms  and find out why it is efficient. We also discuss
the essence of algorithms and its link to self-organizing systems, and finally
we propose some important topics for further research.
\end{abstract}

{\bf Citation details:}
X. S. Yang and S. Deb, Cuckoo search: recent advances and applications,
{\it Neural Computing and Applications}, vol. 24, No. 1, pp. 169-174 (2014).

\section{Introduction}

Optimization concerns many disciplines with a wide range of applications. As time, money and resources are always limited,
optimization is ever-increasingly more important. For example, energy-saving designs and green solutions to many industrial problems
require a paradigm shift in thinking and design practice. On the other hand, people want smart and intelligent products,
and computational intelligence has emerged as a promising area with potentially wide-ranging impact.
Nowadays machine learning often uses optimization algorithms to enhance its learning performance,
while optimization also borrow ideas from machine learning such as statistical learning theory and neural networks.
In this article, we will focus on the introduction of cuckoo search – a powerful, nature-inspired metaheuristic
algorithm for optimization and computational intelligence.

In almost all applications in engineering and industry, we are always trying to optimize something -- whether to
minimize the cost and energy consumption, or to maximize the profit, output, performance and efficiency \cite{Yang,KozielYang,YangAPSO}.
The optimal use of available resources of any sort requires a paradigm shift in scientific thinking,
this is because most real-world applications have far more complicated factors and parameters to affect how the system behaves.
For any optimization problem, the integrated components of the optimization process are the optimization algorithm,
an efficient numerical simulator and a realistic-representation of the physical processes we wish to model and optimize.
This is often a time-consuming process, and in many cases, the computational costs are usually very high. Once we have a good model, the overall computation costs are determined by the optimization algorithms used for search and the numerical solver used for simulation.

Search algorithms are the tools and techniques for achieving optimality of the problem of interest.
This search for optimality is complicated further by the fact that uncertainty almost always presents in the real-world systems. Therefore, we seek not only the optimal design but also robust design in engineering and industry. Optimal solutions, which are not robust enough, are not practical in reality. Suboptimal solutions or good robust solutions are often the choice in such cases.
Optimization problems can be formulated in many ways. For example, the commonly used method of least-squares is a special case of maximum-likelihood formulations. By far the most widely formulation is to write a nonlinear optimization problem as
\be \textrm{minimize } \;\; f_i(\x), (i=1,2,..., M), \ee
subject to the constraints
\be h_j(\x)=0, \quad (j=1,2,..., J), \ee
\be g_k(\x) \le 0, \quad (k=1,2,..., K), \ee
where $f_i, h_j$ and $g_k$ are in general nonlinear functions. Here the design vector or
design variables $\x=(x_1, x_2, ..., x_d)$ can be continuous,
discrete or mixed in $d$-dimensional space. The functions $f_i$
are called objective or cost functions, and when $M>1$, the optimization is multiobjective or multicriteria \cite{Yang}. It is possible to combine different objectives into a single objective, though multiobjective optimization can give far more information and insight into the problem. It is worth pointing out that here we write the problem as a minimization problem, it can also be written as a maximization by simply replacing $f_i(\x)$  by $-f_i(\x)$.

When all functions are nonlinear, we are dealing with nonlinear constrained problems. In some special cases when $f_i, h_j$ and $g_k$ are linear, the problem becomes linear, and we can use the widely linear programming techniques such as the simplex method. When some design variables can only take discrete values (often integers), while other variables are real continuous, the problem is of mixed type, which is often difficult to solve, especially for large-scale optimization problems.

\section{The Essence of an Optimization Algorithm}

\subsection{Optimization Algorithm}

Optimization is a process of searching for the optimal solutions to
a particular problem of interest, and this search process can be carried out
using multiple agents which essentially form a system of evolving agents.
This system can evolve by iterations according to a set of rules or mathematical equations.
Consequently, such a system will show some emergent characteristics, leading
to self-organizing states which correspond to some optima in the search
space. Once the self-organized states are reached, we say the system converges.
Therefore, design of an efficient optimization algorithm is equivalent to
mimicking the evolution of a self-organizing system \cite{Ashby,Keller}.

\subsection{The Essence of an Algorithm}

Mathematically speaking, an algorithm is a procedure to generate outputs for
a given set of inputs. From the optimization point of view, an optimization algorithm
generates a new solution $\x^{t+1}$ to a given problem from a know solution
$\x^t$ at iteration or time $t$. That is
\be \x^{t+1}=\ff{A}(\x^t, \ff{p}(t)), \ee
where $\ff{A}$ is a nonlinear mapping from a given solution $d$-dimensional
vector $\x^t$ to a new solution vector $\x^{t+1}$. The algorithm $\ff{A}$ has
$k$ algorithm-dependent parameters $\ff{p}(t)=(p_1, ..., p_k)$ which can time-dependent
and can thus be tuned.

To find the optimal solution $\x_*$ to a given optimization problem $S$ with an often
infinitely number of states is to select some desired states $\phi$ from all states
$\psi$, according to some predefined criterion $D$. We have
\be S(\psi) \stackrel{\ff{A}(t)}{\longrightarrow} S(\phi(x_*)), \label{algo100} \ee
where the final converged state $\phi$ corresponds to an optimal solution $\x_*$ of the
problem of interest. The selection of the system states in the design space is carried out
by running the optimization algorithm $\ff{A}$. The behavior of the algorithm is
controlled by the parameters $\ff{p}$, the initial solution $\x^{t=0}$ and the
stopping criterion $D$. We can view the combined $S+\ff{A}(t)$ as a complex system
with a self-organizing capability.

The change of states or solutions of the problem of interest is through the algorithm $\ff{A}$.
In many classical algorithms such as hill-climbing, gradient information
of the problem $S$ is used so as to select states, say, the minimum value of the landscape,
and the stopping criterion can be a given tolerance or accuracy, or zero gradient, etc.

An algorithm can act like a tool to tune a complex system.
If an algorithm does not use any state information of the problem, then the algorithm
is more likely to be versatile to deal with many types of problem. However,
such black-box approaches can also imply that the algorithm may  not be efficient as
it could be for a given type of problem. For example, if the optimization problem
is convex, algorithms that use such convexity information will be more efficient
than the ones that do not use such information. In order to be efficient to
select states/solutions efficiently, the information from the search process
should be used to enhance the search process. In many case, such information
is often fed into the selection mechanism of an algorithm. By far the most widely
used selection mechanism to select or keep the best solution found so far. That is,
the simplest form of `survival of the fitness'.

From the schematic representation (\ref{algo100}) of the optimization process, we can see that
the performance of an algorithm may also depend on the type of problem $S$ it solves.
On the other hand, the final, global optimality is achievable or not (within a given
number of iterations) will also depend on the algorithm used. This may be another
way of stating the so-called no-free-lunch theorems.

Optimization algorithms can very diverse with several dozen widely used algorithms.
The main characteristics of different algorithms will only depend on the
actual, nonlinear, often implicit form of $\ff{A}(t)$ and its parameters $\ff{p}(t)$.

\subsection{Efficiency of an Algorithm}

An efficient optimization algorithm is very important to ensure the optimal solutions are reachable. The essence of an algorithm
is a search or optimization process implemented correctly so as to carry out the desired search (though not necessarily efficiently). It can be integrated and linked with other modelling components. There are many optimization algorithms in the literature and no single algorithm is suitable for all problems \cite{Wolpert}.

Algorithms can be classified as deterministic or stochastic. If an algorithm works in a mechanically deterministic manner without any random nature, it is called deterministic. For such an algorithm, it will reach the same final solution if we start with the same initial point. Hill-climbing and downhill simplex are good examples of deterministic algorithms. On the other hand, if there is some randomness in the algorithm, the algorithm will usually reach a different point every time we run the algorithm, even though we start with the same initial point. Genetic algorithms and hill-climbing with a random restart are good examples of stochastic algorithms.

Analyzing current metaheuristic algorithms in more detail, we can single out the type of randomness that a particular algorithm is employing. For example, the simplest and yet often very efficient method is to introduce a random starting point for a deterministic algorithm. The well-known hill-climbing with random restart is a good example. This simple strategy is both efficient in most cases and easy to implement in practice. A more elaborate way to introduce randomness to an algorithm is to use randomness inside different components of an algorithm, and in this case, we often call such  algorithms  heuristic or more often metaheuristic \cite{Yang,YangFA}.

Metaheuristic algorithms are often nature-inspired, and they are now among the most widely used algorithms for optimization. They have many advantages over conventional algorithms \cite{Yang,Gandomi}. Metaheuristic algorithms are very diverse, including genetic algorithms, simulated annealing, differential evolution, ant and bee algorithms, bat algorithm, particle swarm optimization, harmony search, firefly algorithm, cuckoo search and others
\cite{Kennedy,YangFA,YangBA2012,Gandomi2}.
Here we will introduce cuckoo search in great detail.

\section{Cuckoo Search and Analysis}

\subsection{Cuckoo Search}
Cuckoo search (CS) is one of the latest nature-inspired metaheuristic algorithms, developed in 2009 by Xin-She Yang and Suash Deb
\cite{YangDeb,YangDeb2010,YangDeb2012}. CS is based on the brood parasitism of some cuckoo species. In addition, this algorithm is enhanced by the so-called Lévy flights \cite{Pav}, rather than by simple isotropic random walks. Recent studies show that CS is potentially far more efficient than PSO and genetic algorithms \cite{YangDeb}.
Cuckoo are fascinating birds, not only because of the beautiful sounds they can make, but also because of their aggressive reproduction strategy. Some species such as the  {\it ani} and
{\it Guira} cuckoos lay their eggs in communal nests, though they may remove others' eggs to increase the hatching probability of their own eggs. Quite a number of species engage the obligate brood parasitism by laying their eggs in the nests of other host birds (often other species).

For simplicity in describing the standard Cuckoo Search, we now use the following three idealized rules:
\begin{itemize}

\item Each cuckoo lays one egg at a time, and dumps it in a randomly chosen nest;
\item The best nests with high-quality eggs will be carried over to the next generations;
\item The number of available host nests is fixed, and the egg laid by a cuckoo is discovered by the host bird with a probability  . In this case, the host bird can either get rid of the egg, or simply abandon the nest and build a completely new nest.
\end{itemize}

As a further approximation, this last assumption can be approximated by a fraction $p_a$ of the $n$ host nests are replaced by new nests (with new random solutions).
For a maximization problem, the quality or fitness of a solution can simply be proportional to the value of the objective function. Other forms of fitness can be defined in a similar way to the fitness function in genetic algorithms.

For the implementation point of view, we can use the following simple representations that each egg in a nest represents a solution, and each cuckoo can lay only one egg (thus representing one solution), the aim is to use the new and potentially better solutions (cuckoos) to replace a not-so-good solution in the nests. Obviously, this algorithm can be extended to the more complicated case where each nest has multiple eggs representing a set of solutions. For this present introduction, we will use the simplest approach where each nest has only a single egg. In this case, there is no distinction between egg, nest or cuckoo, as each nest corresponds to one egg which also represents one cuckoo.

This algorithm uses a balanced combination of a local random walk and the global explorative random walk,
controlled by
a switching parameter $p_a$. The local random walk can be written as
\be \x_i^{t+1}=\x_i^t +\alpha s \otimes H(p_a-\epsilon) \otimes (\x_j^t-\x_k^t), \ee
where $\x_j^t$ and $\x_k^t$ are two different solutions selected randomly by random permutation,
$H(u)$ is a Heaviside function, $\epsilon$ is a random number drawn from a uniform distribution, and
$s$ is the step size. On the other hand, the global random walk is carried out by using L\'evy flights
\be \x_i^{t+1}=\x_i^t+\a L(s,\lam), \ee
where \be L(s,\lam)=\frac{\lam \Gamma(\lam) \sin (\pi \lam/2)}{\pi}
\frac{1}{s^{1+\lam}}, \quad (s \gg s_0>0). \ee
Here $\a>0$ is the step size scaling factor, which should be related to the scales of the problem of
interest. In most cases,  we can use $\a=O(L/10)$, where $L$ is the characteristic
scale of the problem of interest, while in some cases  $\a=O(L/100)$ can be more effective
 and avoid flying to far. The above equation is essentially the
stochastic equation for a random walk. In general, a random walk is a Markov chain
whose next status/location only depends on the current location (the first term in the above equation)
and the transition probability (the second term).  However,
a substantial fraction of the new solutions should be generated by far field randomization
and their locations should be far enough from the current best solution; this
will make sure that the system will not be trapped in a local optimum \cite{YangDeb,YangDeb2010}.

The literature on cuckoo search is expanding rapidly.
There has been a lot of attention and recent studies using cuckoo search
with a diverse range of applications \cite{Chifu,Choud,Dhiv,Dhiv2,Durgun,Gandomi,Kaveh,Yildiz}.
Walton et al. improved the algorithm by formulating a modified cuckoo
search algorithm \cite{Walton},
while Yang and Deb extended it to multiobjective optimization problems \cite{YangDeb2012}.

\subsection{Why Cuckoo Search is so Efficient?}

Theoretical studies of particle swarm optimization have suggested that
PSO can converge quickly to the current best solution, but not necessarily
the global best solutions \cite{Clerc,Jiang,Wang}. In fact, some analyses suggest
that PSO updating equations do not satisfy the global convergence conditions,
and thus there is no guarantee for global  convergence. On the other hand, it has proved that
cuckoo search satisfy the global convergence requirements and thus has guaranteed global
convergence properties \cite{Wang}. This implies that for multimodal optimization,
PSO may converge prematurely to a local optimum, while cuckoo search can usually
converge to the global optimality.

Furthermore, cuckoo search has two search capabilities: local search and global search,
controlled by a switching/discovery probability. As mentioned in Section 3.1,
the local search is very intensive with about 1/4 of the search time (for
$p_a=0.25$), while global search takes about 3/4 of the total search time.
This allows that the search space can be explored more efficiently on the
global scale, and consequently the global optimality can be found with a higher
probability.

A further advantage of cuckoo search is that its global search uses L\'evy
flights or process, rather than standard random walks. As L\'evy flights have infinite
mean and variance, CS can explore the search space more efficiently than
algorithms by standard Gaussian process. This advantage, combined with
both local and search capabilities and guaranteed global convergence, makes
cuckoo search very efficiently.   Indeed, various studies and applications have demonstrated
that cuckoo search is very efficient \cite{YangDeb2010,Gandomi,Gandomi2,Walton,Civi,Sriv}.

\section{Applications}

Cuckoo search has been applied in many areas of optimization and computational intelligene with promising efficiency. For example, in the engineering design applications, cuckoo search has superior performance over other algorithms for a range of continuous optimization problems such as spring design and welded beam design problems \cite{YangDeb2010,Gandomi,Gandomi2}.

In addition, a modifed cuckoo search by Walton et al. \cite{Walton}
has demonstrated to be very efficient for solving nonlinear problems such as mesh generation. Vazquez \cite{Vazq} used cuckoo search to train spiking neural network models, while Chifu et al. \cite{Chifu} optimized semantic web service composition processes using cuckoo search. Furthermore, Kumar and Chakarverty \cite{Kumar} achieved optimal design for reliable embedded system using cuckoo search, and Kaveh and Bakhshpoori \cite{Kaveh}
used CS to successfully design steel frames. Yildiz \cite{Yildiz}
has used CS to select optimal machine parameters in milling operation with enhanced results, and while Zheng and Zhou \cite{Zheng} provided
a variant of cuckoo search using Gaussian process.

On the other hand, a discrete cuckoo search algorithm has been proposed by Tein and Ramli \cite{Tein} to solve nurse scheduling problems. Cuckoo search has also been used to generate independent paths for software testing and test data generation \cite{Sriv,Perumal,Choud}.
In the context of data fussion and wireless sensor network, cuckoo search has been shown to be very efficient \cite{Dhiv,Dhiv2}.
Furthermore, a variant of cuckoo search in combination with quantum-based approach has been developed to solve Knapsack problems efficiently \cite{Layeb}.
From the algorithm analysis point of view, a conceptural comparison of cuckoo search with particle swarm optimization (PSO), differential evolution (DE), artificial bee colony (ABC) by Civicioglu and Desdo \cite{Civi} suggested
that cuckoo search and differential evoluton algorithms provide more robust results than PSO and ABC. Gandomi et al. \cite{Gandomi} provided a more extensive comparison study for solving various sets of structural optimization problems and concluded that cuckoo search obtained better results than other algorithms such as PSO and gentic algorithms (GA). Speed \cite{Speed} modified the L\'evy cuckoo search and shown that CS can deal with very large-scale problems. Among the diverse applications, an interesting performance enhancement has been obtained by using cuckoo search to train neural networks as shown by Valian et al. \cite{Valian} and reliability optimization problems \cite{Valian2}.

For complex phase equilibrium applications, Bhargava et al. \cite{Bhar} have shown that cuckoo search offers a reliable method for solving thermodynamic calculations. At the same time,
Bulatovi\'c et al. \cite{Bulat} have solved a six-bar double dwell linkage problem using cuckoo search, and Moravej and Akhlaghi \cite{Mora} have solved DG allocation problem in
distribution networks with good convergence rate and performance. Taweewat and Wutiwiwatchi have combined cuckoo search and supervised neural network to estimate musical pitch with reduced size and higher accuracy \cite{Tawee}.

As a further extension, Yang and Deb \cite{YangDeb2012} produced a multiobjective cuckoo search (MOCS)
for design engineering appplications. For multiobjective scheduling problems, a significant progress was made by Chandrasekaran and Simon \cite{Chand} using cuckoo search algorithm, which demonstrated the superiority of their proposed methodology.

Recent studies have demonstrated that cuckoo search can performance significantly better than other algorithms in many applications \cite{Gandomi,Nogh,Zheng,Yildiz}.

\section{Discussion and Concluding Remarks}

Swarm intelligence based algorithms such as cuckoo search and particle swarm optimization are very
efficient in solving a wide range of nonlinear optimization problems, and thus have diverse
applications in sciences and engineering. Some algorithms (e.g., cuckoo search)
can have very good global convergence. However, there are still some challenging issues
that need to be resolved in future studies.

One key issue is that there is a significant gap between theory and practice.
Most metaheuristic algorithms have good applications in practice, but  mathematical
analysis of these algorithms lacks far behind. In fact, apart from a few limited
results about the convergence and stability about algorithms such as particle swarm,
genetic algorithms and cuckoo search, many algorithms do not have theoretical analysis.
Therefore, we may know they can work well in practice, we hardly understand why it works
and how to improve them with a good understanding of their working mechanisms.

Another important issue is that all metaheuristic algorithms have algorithm-dependent parameters,
and the actual values/setting of these parameters will largely influence the performance
of an algorithm. Therefore, the proper parameter-tuning itself
becomes an optimization problem. In fact, parameter-tuning is an important area of research \cite{Eiben},
which deserves more research attention.

In addition, even with very good applications of many algorithms, most of these applications
concern the cases with the number of design variables less than a few hundreds.
It would be more beneficial to real-world applications if the number of variables
can increase to several thousands or even to the order of millions.

All these challenging issues may motivate more research in the near future.
There is no doubt more applications of cuckoo search will be seen in the expanding
literature in the coming years.


\begin{thebibliography}{100}
\bibitem{Ashby}
Ashby, W. R. (1962). Princinples of the self-organizing sysem, in: {\it Pricinples of Self-Organization: Transactions of the
University of Illinois Symposium} (Eds H. Von Foerster and G. W. Zopf, Jr.), Pergamon Press, London, UK. pp. 255-278.


\bibitem{Bhar}
Bhargava, V., Fateen, S. E. K, Bonilla-Petriciolet, A., (2013).
Cuckoo search: a new nature-inspired optimization method
for phase equilibrium calculations, {\it Fluid Phase Equilibria},
{\bf 337}, pp. 191-200.

\bibitem{Bulat}
Bulatovi\'c, R. R., Bordevi\'c, S. R., Dordevi\'c, V. S., (2013).
Cuckoo search algorithm: a metaheuristic approach to solving the
problem of optimum synthesis of a six-bar double dwell linkage,
{\it Mechanism and Machine Theory}, {\bf 61}, pp. 1-13.


\bibitem{Chand}
Chandrasekaran, K., and Simon, S. P., (2012).
Multi-objective scheduling problem: hybrid appraoch using fuzzy assisted
cuckoo search algorithm, {\it Swarm and Evolutionary Computation},
{\bf 5}(1), pp. 1-16.

\bibitem{Chifu}
Chifu, V. R., Pop, C. B.,  Salomie, I., Suia, D. S., and Niculici, A. N., (2012).  Optimizing the semantic web service composition process using cuckoo search, in: Intelligent Distributed Computing V, Studies in Computational Intelligence, {\bf 382}, pp. 93-102.


\bibitem{Choud} Choudhary, K., and Purohit, (2011). G. N., A new testing approach using cuckoo search to achieve multi-objective genetic algorithm, {\it J. of Computing}, {\bf 3}, No. 4, 117-119.

\bibitem{Clerc}
Clerc, M., and Kennedy, J., (2002). The particle swarm --- explosion, stability, and convergence in a multidimensional complex space,
{\it IEEE Trans. Evolutionary Computation}, {\bf 6} (1), pp. 58-73.

\bibitem{Civi} Civicioglu, P., and Besdok, E., (2011).  A conception comparison of the cuckoo search, particle swarm optimization, differential evolution and artificial bee colony algorithms, {\it Artificial Intelligence Review},  doi:10.1007/s10462-011-92760, 6 July (2011).

\bibitem{Dhiv} Dhivya, M., Sundarambal, M., Anand, L. N., (2011). Energy efficient computation of data fusion in wireless sensor networks using cuckoo based particle approach (CBPA), {\it Int. J. of Communications, Network and System Sciences}, {\bf 4}, No. 4, pp. 249-255.

\bibitem{Dhiv2} Dhivya, M. and Sundarambal, M., (2011). Cuckoo search for data gathering in wireless sensor networks, {\it Int. J. Mobile Communications}, {\bf 9}, pp.642-656.


\bibitem{Durgun}
Durgun, I, Yildiz, A.R. (2012). Structural design optimization of vehicle components using cuckoo search algorithm, {\it Materials Testing}, {\bf 3},  185-188.

\bibitem{Eiben}
Eiben A. E. and Smit S. K., (2011).
Parameter tuning for configuring and analyzing evolutionary algorithms, Swarm and Evolutionary Computation, {\bf 1}, pp. 19-31.

\bibitem{Gandomi} Gandomi, A.H., Yang, X.S. and Alavi, A.H. (2013). Cuckoo search algorithm:
a meteheuristic
approach to solve structural optimization problems, {\it Engineering with Computers},
{\bf 29}(1), pp. 17-35 (2013). DOI: 10.1007/s00366-011-0241-y

\bibitem{Gandomi2} Gandomi, A.H., Yang, X.S., Talatahari, S., Deb, S. (2012).
Coupled eagle strategy and differential evolution for unconstrained and constrained global optimization,
{\it Computers \& Mathematics with Applications}, {\bf 63}(1), 191-200.

\bibitem{Jiang}
Jiang, M., Luo, Y. P., and Yang S. Y., (2007). Stochastic convergence
analysis and parameter selection of the standard particle swarm optimization
algorithm, {\it Information Processing Letters}, {\bf 102}, pp. 8-16.


\bibitem{Kaveh}  Kaveh,  A., Bakhshpoori, T., (2011). Optimum design of steel frames using cuckoo search algorithm with Levy flights, Structural Design of Tall and Special Buildings, Vol. 21, online first 28 Nov 2011. http://onlinelibrary.wiley.com/doi/10.1002/tal.754/abstract

\bibitem{Keller}
Keller, E. F. (2009). Organisms, machines, and thunderstorms: a history of self-organization, part two. Complexity,
emergenece, and stable attractors, {\it Historical Studies in the Natural Sciences}, {\bf 39}(1), 1-31.


\bibitem{Kennedy} Kennedy, J.  and  Eberhart, R.C. (1995). Particle swarm optimization,
in: \textit{Proc. of IEEE International Conference on Neural Networks},
Piscataway, NJ, pp. 1942--1948.

\bibitem{KozielYang} Koziel, S. and Yang, X. S., (2011). {\it Computational Optimization, Methods and Algorithms}, Springer, Germany.

\bibitem{Kumar} Kumar A., and Chakarverty, S., (2011). Design optimization for reliable embedded system using Cuckoo Search, in: {\it Proc. of 3rd Int. Conference on Electronics Computer Technology (ICECT2011)}, pp. 564-268.

\bibitem{Layeb}
Layeb, A., (2011). A novel quantum-inspired cuckoo search for Knapsack problems, {\it Int. J. Bio-inspired Computation}, {\bf  3}, No. 5, pp. 297-305.


\bibitem{Mora}
Moravej, Z., Akhlaghi, A., (2013). A novel approach based on cuckoo search for DG allocation
in distribution network, {\it Electrical Power and Energy Systems},
{\bf 44}, pp. 672-679.


\bibitem{Nogh}
Noghrehabadi, A., Ghalambaz, M., and Vosough, A., (2011).  A hybrid power series -- Cuckoo search optimization algorithm to electrostatic deflection of micro fixed-fixed actuators, {\it Int. J. Multidisciplinary Science and Engineering}, {\bf 2}, No. 4, pp. 22-26.



\bibitem{Pav} Pavlyukevich I. (2007). L\'evy flights, non-local search and simulated annealing,
\textit{J. Computational Physics}, {\bf 226}, 1830--1844.

\bibitem{Perumal}
Perumal, K., Ungati, J.M., Kumar, G., Jain, N., Gaurav, R., and Srivastava, P. R., (2011). Test data generation: a hybrid approach using cuckoo and tabu search, Swarm, Evolutionary, and Memetic Computing (SEMCCO2011), Lecture Notes in Computer Sciences, Vol. 7077, pp. 46-54.

\bibitem{Ren}
Ren, Z. H., Wang, J., and Gao Y. L., (2011).
The global convergence analysis
of particle swarm optimization algorithm based on Markov chain,
{\it Control Theory and Applications} (in Chinese), {\bf 28}(4), pp. 462-466.


\bibitem{Speed}
Speed, E. R., (2010). Evolving a Mario agent using cuckoo search and softmax heuristics, Games Innovations Conference (ICE-GIC), pp. 1-7.

\bibitem{Sriv} Srivastava, P. R., Chis, M., Deb, S., and Yang, X. S., (2012). An efficient optimization algorithm for structural software testing, {\it Int. J. Artificial Intelligence}, {\bf  9} (S12), pp. 68-77.

\bibitem{Tawee}
Taweewat, P. and Wutiwiwatchai,C., 2013.
Musical pitch estimation using a supervised single hidden
layer feed-forward neural network, {\it Expert Systems with Applications},
{\bf 40}, pp. 575-589.

\bibitem{Tein}  Tein, L. H. and Ramli,  (2010). R., Recent advancements of nurse scheduling models and a potential path, in: Proc. 6th IMT-GT Conference on Mathematics, Statistics and its Applications (ICMSA 2010), pp. 395-409.

\bibitem{Valian}
Valian, E., Mohanna, S., and Tavakoli, S., (2011). Improved cuckoo search algorithm for feedforward neural network training,
{\it Int. J. Articial Intelligence and Applications},
{\bf 2}, No. 3, 36-43(2011).

\bibitem{Valian2}
Valian, E., Tavakoli, S., Mohanna, S., and Haghi, A., (2013). Improved cuckoo search for
reliability optimization problems, {\it Computer \& Industrial Engineering},
{\bf 64}, pp. 459-468.


\bibitem{Vazq}
Vazquez, R. A., (2011). Training spiking neural models using cuckoo search algorithm,
2011 IEEE Congress on Eovlutionary Computation (CEC'11), pp.679-686.


\bibitem{Walton} Walton, S., Hassan, O., Morgan, K., and Brown, M.R. (2011).
Modified cuckoo search: a new gradient free optimization algorithm.
{\it Chaos, Solitons \& Fractals}, {\bf 44}(9), 710-718.

\bibitem{Wang}
Wang, F., He, X.-S., Wang, Y., and Yang, S. M., (2012).
Markov model and convergence analysis based on cuckoo search algorithm,
{\it Computer Engineering}, {\bf 38}(11), pp. 180-185.


\bibitem{Wolpert}  Wolpert, D. H. and Macready, W. G.,  No free lunch theorems for optimization, {\it IEEE Trans. Evolutionary Computation},  {\bf 1}, 67--82 (1997)


\bibitem{Yang}
 Yang X. S., (2010). {\it Engineering Optimisation: An Introduction with
 Metaheuristic Applications}, John Wiley and Sons.

\bibitem{YangFA} Yang X. S., (2009).
Firefly algorithms for multimodal optimization, in: Stochastic Algorithms: Foundations and Applications, SAGA 2009, Lecture Notes in Computer Sciences,
Vol. {\bf 5792}, 169-178.

\bibitem{YangFA2} Yang X.-S., (2010).
Firefly algorithm, stochastic test functions and design optimisation, {\it Int. J. Bio-inspired Computation},
{\bf 2}(2),  78-84.

\bibitem{YangAPSO}
Yang, X. S., Deb, S., and Fong, S., (2011). Accelerated particle swarm
optimization and support vector machine for business optimization and
applications, in: Networked Digital Technologies 2011, {\it Communications in
Computer and Information Science}, {\bf 136}, pp. 53--66.


\bibitem{YangBA2012}
Yang, X. S., Gandomi, A. H., (2012). Bat algorithm: a novel approach for global engineering optimization, Engineering Computations, 29(5), 1-18.


\bibitem{YangDeb}
Yang, X. S. and Deb, S., (2009). Cuckoo search via L\'evy flights,
{\it Proceeings of World Congress on Nature \& Biologically Inspired
Computing} (NaBIC 2009), IEEE Publications, USA, pp. 210-214.

\bibitem{YangDeb2010}
Yang X. S. and Deb S., (2010). Engineering optimization by cuckoo search,
{\it Int. J. Math. Modelling Num. Opt.}, {\bf 1} (4), 330-343 (2010).

\bibitem{YangDeb2012} Yang, X.S. and Deb, S. (2012).  Multiobjective cuckoo search for design optimization,
{\it Computers and Operations Research}, accepted October (2011). doi:10.1016/j.cor.2011.09.026


\bibitem{Yildiz} Yildiz, A. R., (2012).  Cuckoo search algorithm for the selection of optimal machine parameters in milling operations, {\it Int. J. Adv. Manuf. Technol.}, (2012). doi:10.1007/s00170-012-4013-7

\bibitem{Zheng}
Zheng, H. Q., and Y. Zhou, Y. (2012). A novel cuckoo search optimization algorithm based on Gauss distribution, {\it J. Computational Information Systems}, {\bf 8}, pp.4193-4200.


\end{thebibliography}
\end{document}